\documentclass[10pt]{article}
\usepackage{amsfonts}
\def\frak{\fontencoding{U}\fontfamily{euf}\selectfont}
\begin{document}
\topmargin= -.2in \baselineskip=16pt

\title {$L$-functions
of Symmetric Products of the Kloosterman Sheaf over {\bf Z}
}

\author {Lei Fu \thanks{The research of Lei Fu is supported by
the NSFC (10525107).} \\{\small Chern Institute of Mathematics and
LPMC, Nankai University, Tianjin, P. R.
China}\\
{\small leifu@nankai.edu.cn}\\{}\\
Daqing Wan
\\
{\small Department of Mathematics, University of California, Irvine,
CA
92697}\\
{\small dwan@math.uci.edu}}
\date{}
\maketitle

\begin{abstract}
The classical $n$-variable Kloosterman sums over the finite field
${\bf F}_p$ give rise to a lisse $\overline {\bf Q}_l$-sheaf ${\rm
Kl}_{n+1}$ on ${\bf G}_{m, {\bf F}_p}={\bf P}^1_{{\bf
F}_p}-\{0,\infty\}$, which we call the Kloosterman sheaf. Let
$L_p({\bf G}_{m,{\bf F}_p}, {\rm Sym}^k{\rm Kl}_{n+1}, s)$ be the
$L$-function of the  $k$-fold symmetric product of ${\rm Kl}_{n+1}$.
We construct an explicit virtual scheme $X$ of finite type  over
${\rm Spec}\, {\bf Z}$ such that the $p$-Euler factor of the zeta
function of $X$ coincides with $L_p({\bf G}_{m,{\bf F}_p}, {\rm
Sym}^k{\rm Kl}_{n+1}, s)$. We also prove similar results for
$\otimes^k {\rm Kl}_{n+1}$ and $\bigwedge^k {\rm Kl}_{n+1}$.
\end{abstract}
\bigskip
\centerline {\bf 0. Introduction}

\bigskip
\bigskip
For each prime number $p$, let ${\bf F}_p$ be a finite field with
$p$ elements. Fix an algebraic closure $\overline {\bf F}_p$ of
${\bf F}_p$. For any power $q$ of $p$, let ${\bf F}_q$ be the
subfield of $\overline {\bf F}_p$ with $q$ elements. Let $l$ be a
prime number distinct from $p$. Fix a nontrivial additive character
$\psi: {\bf F}_p\to \overline {\bf Q}_l^\ast$. Thus, $\psi(1)$ is a
primitive $p$-th root of unity, which is denoted by $\zeta_p$. For
any nonzero $x\in {\bf F}_q$, we define the $n$-variable Kloosterman
sum by
$${\rm Kl}_{n+1}({\bf F}_q, x)=\sum\limits_{x_1,\ldots, x_{n+1}
\in  {\bf F}_q^\ast , x_1\cdots x_{n+1}=x } \psi({\rm Tr}_{{\bf
F}_q/{\bf F}_p}(x_1+\cdots + x_{n+1}))\in {\bf Z}[\zeta_p].$$ In
[SGA $4\frac{1}{2}$] [Sommes trig.] \S 7, Deligne constructs a lisse
$\overline {\bf Q}_l$-sheaf ${\rm Kl}_{n+1}$ on ${\bf G}_{m, {\bf
F}_p}={\bf P}^1_{{\bf F}_p}-\{0,\infty\}$ such that for any $x\in
{\bf G}_m({\bf F}_q)={\bf F}_q^\ast$, we have
$${\rm Tr}(F_x, {\rm Kl}_{n+1,\bar x})=(-1)^n {\rm Kl}_{n+1}({\bf F}_q,
x),$$ where $F_x$ is the geometric Frobenius element at the point
$x$. For any natural number $k$, consider the $L$-function
$$L_p({\bf G}_{m, {\bf F}_p}, {\rm Sym}^k {\rm Kl}_{n+1}, s)
=\prod_{x\in |{\bf G}_{m, {\bf F}_p}|} {\rm det}(1-F_x p^{-s\, {\rm
deg}(x)}, ({\rm Sym}^k {\rm Kl}_{n+1})_{\bar x})^{-1}$$ of the
$k$-fold symmetric product ${\rm Sym}^k{\rm Kl}_{n+1}$ of ${\rm
Kl}_{n+1}$, where $|{\bf G}_{m,{\bf F}_p}|$ is the set of Zariski
closed points in ${\bf G}_{m,{\bf F}_p}$. This $L$-function in $s$
has two parameters $k$ and $p$. It was first studied by Robba [Ro]
in the case $n=1$ via $p$-adic methods. More recently,  its basic
properties and $p$-adic variation as $k$ varies $p$-adicaly have
been studied extensively in connection with Dwork's unit root
conjecture. See [W1], [GK], [FW1] and [FW2].

In this paper, we fix $k$ and study the variation of this
$L$-function as $p$ varies. It was observed in [FW1] Lemma 2.2 that
for each $p$, the $L$-function $L_p({\bf G}_{m, {\bf F}_p}, {\rm
Sym}^k {\rm Kl}_{n+1}, s)$ is a polynomial in $p^{-s}$ with
coefficients in ${\bf Z}$. This naturally leads to the conjecture
that the infinite product
$$\zeta_{k,n}(s):=\prod_p L_p({\bf G}_{m,{\bf F}_p}, {\rm Sym}^k{\rm Kl}_{n+1}, s)$$
is automorphic and thus extends to a meromorphic function in $s\in
{\bf C}$. This is easy to prove if $n=1$ and $k\leq 4$. If $n=1$ and
$k=5$, the series $\zeta_{5,1}(s)$ is essentially the $L$-function
of an elliptic curve with complex multiplication and thus meromprhic
in $s\in {\bf C}$, see [PTV]. If $n=1$ and $k=6$, the modularity of
$\zeta_{6,1}(s)$ follows from [HS] and the references listed there.
In this case, one obtains a rigid Calabi-Yau threefold. In the case
$n=1$ and $k=7$, the series $\zeta_{7,1}(s)$ is conjectured by Evans
[Ev] to be given by the $L$-function associated to an explicit
modular form of weight $3$ and level $525$. With the recent progress
on the modularity problem due to Taylor and Harris, it may be
possible to prove the meromorphic continuation of $\zeta_{k,n}(s)$
for some larger $k$ and $n$.

To prove the meromorphic continuation of $\zeta_{k,n}(s)$, the first step would be
to prove that $\zeta_{k,n}(s)$ is motivic (or geometric) in nature, i.e., it arises
as the zeta function of a motive over ${\rm Spec}{\bf Z}$.
This question was raised in [FW1] and is solved in this paper.
We will construct a virtual $\overline {\bf Q}_l$-sheaf ${\cal G}$ of
geometric origin on ${\rm Spec}\, {\bf Z}$ so that the Euler factor
$$L_p({\rm Spec}\, {\bf Z},{\cal G},s)={\rm det}(1-F_p p^{-s},
{\cal G}_{\bar p})^{-1}$$ of the $L$-function of ${\cal G}$
coincides with $L_p({\bf G}_{m, {\bf F}_p}, {\rm Sym}^k{\rm
Kl}_{n+1}, s)$ for each prime number $p$, where $F_p$ is the
geometric Frobenius element at $p$. We also prove similar results
for $\otimes ^k {\rm Kl}_{n+1}$ and $\bigwedge^k {\rm Kl}_{n+1}$.

To describe our results, we introduce the following schemes over
${\bf Z}$.

\bigskip
\noindent {\bf Definition 0.1}. Denote the homogeneous coordinates
of ${\bf P}^{kn-1}$ by $[x_{ij}]$ $(i=1,\ldots, n,\; j=1,\ldots,
k)$. Let $Y_k$ be the subscheme of ${\bf P}^{kn-1}$ defined by
$$x_{ij}\not=0,$$let $Y_{k0}$ be the subscheme defined by
$$x_{ij}\not=0,\;\sum_{i,j}
x_{ij}=0,$$ let $Z_k$ be the subscheme defined by the conditions
$$x_{ij}\not=0,\; \sum_{j=1}^k \frac{1}{\prod_{i=1}^n x_{ij}}=0,$$
and let $Z_{k0}$ be the subscheme defined by
$$x_{ij}\not=0,\; \sum_{j=1}^k \frac{1}{\prod_{i=1}^n x_{ij}}=0,\;
\sum_{i,j} x_{ij}=0.$$ These are schemes of finite type over ${\bf
Z}$. Let $\hbox {\frak S}_k$ be the group of permutations of the set
$\{1,\ldots, k\}$. It acts on ${\bf P}^{kn-1}$ by permuting the
homogenous coordinates $x_{i1},\ldots, x_{ik}$ for each $i$.
Similarly $\hbox {\frak S}_k$ acts on $Z_{k0}$, $Z_k$, $Y_{k0}$ and
$Y_k$. The notations $Y_k/\hbox{\frak S}_k$, $Y_{k0}/\hbox{\frak
S}_k$, $Z_k/\hbox{\frak S}_k$, and $Z_{k0}/\hbox{\frak S}_k$ denote
the quotient scheme of $Y_k$, $Y_{k0}$, $Z_k$, and $Z_{k0}$ by
$\hbox{\frak S}_k$, respectively. Our main result is the following
theorem.

\bigskip
\noindent {\bf Theorem 0.2.} For a scheme $X$ of finite type over
${\bf Z}$, let $\zeta_X(s)$ denote its zeta function. We have
$$L_p({\bf G}_{m,{\bf F}_p}, \otimes^k{\rm Kl}_{n+1}, s)
=\left(\frac{\zeta_{Z_{k0,{\bf F}_p}}(s-2)\zeta_{Y_{k,{\bf
F}_p}}(s)}{\zeta_{Z_{k,{\bf F}_p}}(s-1)\zeta_{Y_{k0,{\bf
F}_p}}(s-1)}\right)^{(-1)^{kn}},$$
$$
L_p({\bf G}_{m,{\bf F}_p}, {\rm Sym}^k{\rm Kl}_{n+1},
s)=\left(\frac{\zeta_{Z_{k0,{\bf F}_p}/\hbox{\frak
S}_k}(s-2)\zeta_{Y_{k,{\bf F}_p}/\hbox{\frak
S}_k}(s)}{\zeta_{Z_{k,{\bf F}_p}/\hbox{\frak
S}_k}(s-1)\zeta_{Y_{k0,{\bf F}_p}/\hbox{\frak
S}_k}(s-1)}\right)^{(-1)^{kn}}.$$ Thus,
$$\prod_p L_p({\bf G}_{m,{\bf F}_p}, \otimes^k{\rm Kl}_{n+1}, s)=
\left(\frac{\zeta_{Z_{k0}}(s-2)\zeta_{Y_k}(s)}
{\zeta_{Z_{k}}(s-1)\zeta_{Y_{k0}}(s-1)}\right)^{(-1)^{kn}},$$
$$\prod_p L_p({\bf
G}_{m,{\bf F}_p}, {\rm Sym}^k{\rm Kl}_{n+1},
s)=\left(\frac{\zeta_{Z_{k0}/\hbox{\frak
S}_k}(s-2)\zeta_{Y_{k}/\hbox{\frak
S}_k}(s)}{\zeta_{Z_{k}/\hbox{\frak
S}_k}(s-1)\zeta_{Y_{k0}/\hbox{\frak
S}_k}(s-1)}\right)^{(-1)^{kn}}.$$
\bigskip

The above formulas can be simplified significantly. This is done in
\S 4. To prove the above results, we need to relate Kloosterman
sheaves by the $l$-adic Fourier transformation. This is done in \S
1. We prove Theorem 0.2 in \S 2 and \S 3.

\bigskip
\noindent {\bf Remark 0.3.} For any partition $\lambda$ of $k$, let
$S_\lambda({\rm Kl}_{n+1})$ be the Weyl construction applied to
${\rm Kl}_{n+1}$. (Confer [FH] \S 6.1.) The method developed in this
paper can also be used to show that $L_p({\bf G}_{m,{\bf
F}_p},S_\lambda({\rm Kl}_{n+1}), s)$ is the Euler factor at $p$ of
the $L$-function of a virtual $\overline{\bf Q}_l$-sheaf on ${\rm
Spec}\,{\bf Z}$ of geometric origin for each prime number $p$. An
example is given in Theorem 3.2 for the $k$-th exterior product.

\bigskip
\centerline {\bf 1. Kloosterman Sheaves and the Fourier
Transformation}

\bigskip
\bigskip
In this section, we give an inductive construction of Kloosterman
sheaves using the $l$-adic Fourier transformation. We refer the
reader to [L] for the definition and properties of the Fourier
transformation.

The morphism
$${\cal P}: {\bf A}_{{\bf F}_p}^1\to {\bf A}_{{\bf F}_p}^1$$
corresponding to the ${\bf F}_p$-algebra homomorphism
$${\bf F}_p[t]\to {\bf F}_p[t],\; t\mapsto t^p-t$$
is a finite galois \'etale covering space, and it defines an ${\bf
F}_p$-torsor
$$0\to {\bf F}_p\to {\bf A}_{{\bf F}_p}^1\stackrel {\cal P}\to
{\bf A}_{{\bf F}_p}^1\to 0.$$ Pushing-forward this torsor by
$\psi^{-1}:{\bf F}_p\to\overline{\bf Q}_l$, we get a lisse
$\overline{\bf Q}_l$-sheaf ${\cal L}_\psi$ of rank $1$ on ${\bf
A}_{{\bf F}_p}^1$, which we call the Artin-Schreier sheaf. Let $X$
be a scheme over ${\bf F}_p$ and let $f$ be an element in the ring
of global sections $\Gamma(X,{\cal O}_X)$ of the structure sheaf of
$X$. Then $f$ defines an ${\bf F}_p$-morphism $X\to {\bf A}_{{\bf
F}_p}^1$ so that the induced ${\bf F}_p$-algebra homomorphism ${\bf
F}_p[t]\to \Gamma(X,{\cal O}_X)$ maps $t$ to $f$. We often denote
this canonical morphism also by $f$, and denote by ${\cal
L}_{\psi}(f)$ the inverse image of ${\cal L}_\psi$ under this
morphism.

The main result of this section is the following.

\bigskip
\noindent{\bf Proposition 1.1.} Let $i:{\bf G}_{m,{\bf F}_p}\to{\bf
G}_{m,{\bf F}_p}$ be the morphism $x\mapsto \frac{1}{x}$, and let
$j:{\bf G}_{m,{\bf F}_p}\to {\bf A}^1_{{\bf F}_p}$ be the canonical
open immersion. For each integer $n\geq 1$, define ${\rm Kl}_n$
inductively as follows:
\begin{eqnarray*}
{\rm Kl}_1&=&{\cal L}_\psi|_{{\bf G}_{m,{\bf F}_p}},\\
{\rm Kl}_{n+1}&=& ({\cal F}(j_! i^\ast {\rm Kl}_{n}))|_{{\bf G}_{m,
{\bf F}_p}},
\end{eqnarray*}
where ${\cal F}(-)=Rp_{2!}(p_1^\ast(-)\otimes^L {\cal
L}_\psi(tt'))[1]$ denotes the Fourier transformation. Here
$$p_1,p_2:{\bf A}^1_{{\bf F}_p}\times_{{\bf F}_p}{\bf A}^1_{{\bf
F}_p}\to {\bf A}^1_{{\bf F}_p}$$ are the projections, and $tt'$ is
regarded as an element in $$\Gamma({\bf A}^1_{{\bf F}_p}\times_{{\bf
F}_p}{\bf A}^1_{{\bf F}_p},{\cal O}_{{\bf A}^1_{{\bf
F}_p}\times_{{\bf F}_p}{\bf A}^1_{{\bf F}_p}})\cong {\bf
F}_p[t,t'].$$

(i) For any $t\in{\bf G}_m({\bf F}_q)$, we have
$${\rm Tr}(F_t, {\rm Kl}_{n,\bar t})=(-1)^{n-1}\sum\limits_{x_1,\ldots, x_n
\in  {\bf F}_q^\ast,\; x_1\cdots x_n=t } \psi({\rm Tr}_{{\bf
F}_q/{\bf F}_p}(x_1+\cdots + x_n)).$$

(ii) ${\rm Kl}_n$ is a lisse $\overline{\bf Q}_l$-sheaf on ${\bf
G}_{m,{\bf F}_p}$ of rank $n$. It is tame at $0$, and its Swan
conductor at $\infty$ is $1$.

\bigskip
It follows from the proposition that the sheaf ${\rm Kl}_n$ defined
inductively using the Fourier transformation as above coincides with
the Kloosterman sheaf constructed by Deligne.

\bigskip
\noindent {\bf Proof.} We use induction on $n$. When $n=1$, the
assertions are clear. Suppose the assertions hold for ${\rm Kl}_n$.
We have
\begin{eqnarray*}
&&{\rm Tr}(F_t, {\rm Kl}_{n+1,\bar t})\\
&=&{\rm Tr}(F_t, ({\cal F}(j_! i^\ast {\rm Kl}_{n}))_{\bar t})\\
&=& -\sum_{s\in {\bf F}_q}\psi({\rm Tr}_{{\bf F}_q/{\bf
F}_p}(st)){\rm Tr}(F_s,(j_! i^\ast {\rm Kl}_{n})_{\bar s})\\
&=&(-1)^n \sum_{s\in {\bf F}_q^\ast}\psi({\rm Tr}_{{\bf F}_q/{\bf
F}_p}(st))\sum\limits_{x_1,\ldots, x_{n}\in  {\bf
F}_q^\ast,\;x_1\cdots x_n=\frac{1}{s} } \psi({\rm Tr}_{{\bf
F}_q/{\bf F}_p}(x_1+\cdots +
x_n))\\
&=&(-1)^n \sum_{s,x_1,\ldots, x_n\in  {\bf F}_q^\ast,\; x_1\cdots
x_n=\frac{1}{s}} \psi({\rm Tr}_{{\bf F}_q/{\bf F}_p}(x_1+\cdots +
x_n+st))\\
&=&(-1)^n \sum_{x_1,\ldots, x_{n+1}\in  {\bf F}_q^\ast,\; x_1\cdots
x_{n+1}=t } \psi({\rm Tr}_{{\bf F}_q/{\bf F}_p}(x_1+\cdots +
x_{n+1})),
\end{eqnarray*}
where the second equality follows from the definition of the Fourier
transformation, and the third equality follows from the induction
hypothesis. This proves (i) holds for ${\rm Kl}_{n+1}$. Let $\eta_0$
(resp. $\eta_\infty$) be the generic point of the strict
henselization of ${\bf P}_{{\bf F}_p}^1$ at $0$ (resp. $\infty$). By
the induction hypothesis, ${\rm Kl}_n$ is tame at $0$. Hence
$(i^\ast {\rm Kl}_n))|\eta_{\infty}$ is tame. By [L] 2.3.1.3 (i),
${\rm Kl}_{n+1}=({\cal F}(j_! i^\ast {\rm Kl}_{n}))|_{{\bf G}_{m,
{\bf F}_p}}$ is a lisse sheaf on ${\bf G}_{m,{\bf F}_p}$. Moreover,
by [L] 2.5.3.1, ${\cal F}^{(\infty,0')} ((i^\ast {\rm
Kl}_n)|\eta_\infty)$ is tame. It follows that ${\rm Kl}_{n+1}$ is
tame at $0$. By the stationary phase principle [L] 2.3.3.1 (iii), we
have
$${\rm Kl}_{n+1}|\eta_{\infty'}={\cal F}^{(0,\infty')} ((i^\ast
{\rm Kl}_n)|\eta_{0}) \oplus {\cal F}^{(\infty,\infty')} ((i^\ast
{\rm Kl}_n)|\eta_{\infty}).$$  Since $(i^\ast {\rm
Kl}_n))|\eta_{\infty}$ is tame, we have ${\cal F}^{(\infty,\infty')}
((i^\ast {\rm Kl}_n))|\eta_\infty)=0 $ by [L] 2.4.3 (iii) b). By the
induction hypothesis, the Swan conductor of $(i^\ast {\rm
Kl}_n)|\eta_{0}$ is 1 and its rank is $n$. By [L] 2.4.3 (i) b), the
Swan conductor of ${\cal F}^{(0,\infty')} ((i^\ast {\rm
Kl}_n)|\eta_{0})$ is 1, and its rank is $n+1$. Hence the Swan
conductor of ${\rm Kl}_{n+1}$ at $\infty$ is $1$, and the rank of
${\rm Kl}_{n+1}$ is $n+1$. This proves (ii) holds for ${\rm
Kl}_{n+1}$.

\bigskip
\bigskip
\centerline {\bf 2. The $L$-function of $\otimes^k {\rm Kl}_{n+1}$}

\bigskip
\bigskip
Let $$\tilde {\bf A}^{n+1}_{{\bf F}_p}=\{(x,y)\in {\bf A}_{{\bf
F}_p}^{n+1}\times _{{\bf F}_p}{\bf P}_{{\bf F}_p}^n|x \hbox { lies
on the line determined by } y\}$$ be the blowing-up of ${\bf
A}_{{\bf F}_p}^{n+1}$ at the origin, let
$$\pi_1: \tilde {\bf A}_{{\bf F}_p}^{n+1}\to {\bf A}_{{\bf F}_p}^{n+1},\;
\pi_2: \tilde {\bf A}_{{\bf F}_p}^{n+1}\to {\bf P}_{{\bf F}_p}^n$$
be the projections, let
$$H=\{[x_0:\ldots:x_n]\in {\bf P}^n|\sum x_i=0\},$$
and let
$$\kappa: H\to {\bf P}_{{\bf F}_p}^n$$ be the canonical closed immersion.
Consider the morphism $$s: {\bf A}_{{\bf F}_p}^{n+1}\to {\bf
A}_{{\bf F}_p}^1,\; s(x_0,\ldots, x_{n})=x_0+\cdots+ x_{n}.$$ We
have
\begin{eqnarray*}
R\pi_{2!}\pi_1^\ast s^\ast {\cal L}_\psi=\kappa_! \overline {\bf
Q}_l(-1)[-2].
\end{eqnarray*}
This follows from the fact that $\tilde {\bf A}_{{\bf F}_p}^{n+1}$
is a line bundle over ${\bf P}_{{\bf F}_p}^n$, and that for any
point $a=[a_0:\ldots :a_n]$ in ${\bf P}_{{\bf F}_p}^n$, we have
$$R\Gamma_c(\pi_2^{-1}(a)\otimes \overline {\bf F}_p,\pi_1^\ast s^\ast {\cal L}_\psi)\cong
R\Gamma_c({\bf A}^1_{\overline {\bf F}_p},{\cal L}_\psi(t\sum
a_i))=\left\{
\begin{array}{cl}
0 &\hbox { if } \sum a_i \not=0, \\
\overline {\bf Q}_l(-1)[-2] & \hbox { otherwise.}
\end{array}\right.$$

\bigskip
\noindent {\bf Lemma 2.1.} For a subscheme $Z$ of ${\bf P}_{{\bf
F}_p}^n$, let
$$Z_0=Z\cap H,\;\tilde X=\pi_2^{-1}(Z),\; X=\pi_1(\tilde X).$$
We have a natural distinguished triangle
$$R\Gamma_c((X-\{0\})\otimes \overline {\bf F}_p, s^\ast {\cal
L}_\psi)\to R\Gamma_c(Z_0\otimes \overline {\bf F}_p, \overline{\bf
Q}_l(-1)[-2])\to R\Gamma_c(Z\otimes \overline {\bf F}_p,
\overline{\bf Q}_l)\to.$$

\bigskip
\noindent {\bf Proof.} Let $\pi_1':\tilde X\to X$ be the restriction
of $\pi_1$ to $\tilde X$. We have a distinguished triangle
$$R\Gamma_c(\pi_1'^{-1}(X-\{0\})\otimes \overline {\bf F}_p, \pi_1^\ast s^\ast {\cal
L}_\psi)\to R\Gamma_c(\tilde X\otimes \overline {\bf F}_p,
\pi_1^\ast s^\ast {\cal L}_\psi)\to
R\Gamma_c(\pi_1'^{-1}(\{0\})\otimes \overline {\bf F}_p, \pi_1^\ast
s^\ast {\cal L}_\psi)\to.$$ On the other hand, we have
$$X-\{0\}\cong\pi_1'^{-1}(X-\{0\}),\;\pi_1'^{-1}(X)=\tilde X,\;
\pi_1'^{-1}(\{0\})\cong Z,$$ and hence
\begin{eqnarray*}
R\Gamma_c(\pi_1'^{-1}(X-\{0\})\otimes \overline {\bf F}_p,
\pi_1^\ast s^\ast {\cal L}_\psi)&\cong& R\Gamma_c((X-\{0\})\otimes
\overline {\bf F}_p, s^\ast {\cal L}_\psi),\\
R\Gamma_c(\tilde X\otimes \overline {\bf F}_p, \pi_1^\ast s^\ast
{\cal L}_\psi)&\cong& R\Gamma_c(Z\otimes \overline {\bf F}_p,
R\pi_{2!} \pi_1^\ast s^\ast {\cal L}_\psi)\\
&\cong& R\Gamma_c(Z\otimes \overline {\bf F}_p,\kappa_!
\overline{\bf
Q}_l(-1)[-2])\\
&\cong& R\Gamma_c(Z_0\otimes \overline {\bf F}_p, \overline{\bf
Q}_l(-1)[-2]),\\
R\Gamma_c(\pi_1'^{-1}(\{0\})\otimes \overline {\bf F}_p, \pi_1^\ast
s^\ast {\cal L}_\psi)&\cong& R\Gamma_c(Z\otimes \overline {\bf F}_p,
\overline{\bf Q}_l).
\end{eqnarray*}
Our assertion follows.

\bigskip
By Proposition 1.1, we have
$${\cal F}(j_! i^\ast {\rm
Kl}_{n})|_{{\bf G}_{m, {\bf F}_p}}\cong {\rm Kl}_{n+1}.$$ By [L]
1.2.2.7, we have
$${\cal F}(\ast^k (j_! i^\ast {\rm
Kl}_{n}))|_{{\bf G}_{m, {\bf F}_p}}\cong \otimes^k{\rm
Kl}_{n+1}[1-k],\eqno(1)$$ where $\ast^k$ denotes the $k$-fold
convolution product. Let
\begin{eqnarray*}
s_n:{\bf G}_m^n&\to&{\bf A}^1,\\
p_n:{\bf G}_m^n&\to&{\bf G}_m
\end{eqnarray*}
be the morphisms
\begin{eqnarray*}
s_n(x_1,\ldots,x_n)=x_1+\cdots+x_n,\\
p_n(x_1,\ldots,x_n)=x_1\cdots x_n,
\end{eqnarray*}
respectively. By [SGA $4\frac{1}{2}$] [Sommes trig.] \S 7, we have
$${\rm Kl}_n\cong Rp_{n!}s_n^\ast {\cal L}_\psi[n-1].\eqno(2)$$
Denote the coordinates of ${\bf G}_m^{kn}$ by $x_{ij}$ $(i=1,\ldots,
n, \; j=1,\ldots, k)$. Let
\begin{eqnarray*}
s_{kn}:{\bf G}_m^{kn}&\to&{\bf A}^1,\\
f_{kn}:{\bf G}_m^{kn}&\to&{\bf A}^1
\end{eqnarray*}
be the morphisms
\begin{eqnarray*}
s_{kn}((x_{ij}))=\sum_{i,j}x_{ij},\\
f_{kn}((x_{ij}))=\sum_{j=1}^k \frac{1}{\prod_{i=1}^n x_{ij}},
\end{eqnarray*}
respectively. By the K\"{u}nneth formula, the definition of the
convolution product [L] 1.2.2.6, and the isomorphism (2), we have
$$\ast^k (j_! i^\ast {\rm
Kl}_{n})\cong Rf_{kn,!}s_{kn}^\ast {\cal L}_\psi[k(n-1)].$$ Combined
with the isomorphism (1), we get
$$({\cal F}(Rf_{kn,!}s_{kn}^\ast {\cal L}_\psi)[kn-1])|_{{\bf G}_{m, {\bf F}_p}}\cong \otimes^k{\rm
Kl}_{n+1}. \eqno(3)$$ By Grothendieck's formula for $L$-functions,
we have
$$L_p({\bf G}_{m,{\bf F}_p},\otimes
^k {\rm Kl}_{n+1},s)={\rm det}(1-F p^{-s}, R\Gamma_c ({\bf
G}_{m,\overline{\bf F}_p}, \otimes ^k{\rm Kl}_{n+1}))^{-1}.$$ Taking
into account of the isomorphism (3), we get
\begin{eqnarray*}
L_p({\bf G}_{m,{\bf F}_p},\otimes ^k {\rm Kl}_{n+1},s)&=&{\rm
det}(1-F p^{-s}, R\Gamma_c({\bf G}_{m,\overline{\bf F}_p}, ({\cal
F}(Rf_{kn,!}s_{kn}^\ast {\cal L}_\psi)[kn-1])|_{{\bf G}_{m,\overline{\bf F}_p}}))^{-1}\\
&=& \frac{{\rm det}(1-F p^{-s}, R\Gamma_c({\bf A}^1_{\overline{\bf
F}_p}, {\cal F}(Rf_{kn,!}s_{kn}^\ast {\cal
L}_\psi)[kn-1]))^{-1}}{{\rm det}(1-Fp^{-s}, ({\cal
F}(Rf_{kn,!}s_{kn}^\ast {\cal L}_\psi)[kn-1])|_{\bar 0})^{-1}}
\end{eqnarray*}
By the definition of the Fourier transformation, we have $$({\cal
F}(Rf_{kn,!}s_{kn}^\ast {\cal L}_\psi))|_{\bar 0}\cong R\Gamma_c
({\bf A}^1_{\overline{\bf F}_p}, Rf_{kn,!}s_{kn}^\ast {\cal
L}_\psi)[1]\cong R\Gamma_c({\bf G}_{m,\overline{\bf
F}_p}^{kn},s_{kn}^\ast {\cal L}_\psi)[1].$$ Hence
\begin{eqnarray*}
{\rm det}(1-Fp^{-s}, ({\cal F}(Rf_{kn,!}s_{kn}^\ast {\cal
L}_\psi)[kn-1])|_{\bar 0})^{-1}&=&{\rm det}(1-Fp^{-s},
R\Gamma_c({\bf G}^{kn}_{m,\overline{\bf F}_p}, s_{kn}^\ast {\cal
L}_\psi)[kn])^{-1}.
\end{eqnarray*}
By the inversion formula for the Fourier transformation [L] 1.2.2.1,
we have
\begin{eqnarray*}
R\Gamma_c({\bf A}_{\overline{\bf F}_p}^1, {\cal
F}(Rf_{kn,!}s_{kn}^\ast {\cal L}_\psi))
&\cong& ({\cal F}({\cal F}(Rf_{kn,!}s_{kn}^\ast {\cal L}_\psi)))_{\bar 0}[-1] \\
&\cong & (Rf_{kn,!}s_{kn}^\ast {\cal L}_\psi)_{\bar 0}(-1)[-1]\\
&\cong& R\Gamma_c(X_{k,\overline{\bf F}_p},s_{kn}^\ast {\cal
L}_\psi)(-1)[-1],
\end{eqnarray*}
where $X_k$ is the subscheme of ${\bf G}_m^{kn}$ over ${\bf Z}$
defined by the equation
$$\sum_{j=1}^k \frac{1}{\prod_{i=1}^n x_{ij}}=0.$$
Hence
\begin{eqnarray*}
&&{\rm det}(1-F p^{-s}, R\Gamma_c({\bf A}^1_{\overline{\bf F}_p},
{\cal F}(Rf_{kn,!}s_{kn}^\ast {\cal L}_\psi)[kn-1]))^{-1} \\
&=&{\rm det}(1-F p^{-s}, R\Gamma_c(X_{k,\overline{\bf
F}_p},s_{kn}^\ast {\cal L}_\psi)(-1)[kn-2])^{-1}.
\end{eqnarray*}
It follows that
\begin{eqnarray*}
L_p({\bf G}_{m,{\bf F}_p}, \otimes^k{\rm Kl}_{n+1}, s)&=&\frac{{\rm
det}(1-F p^{-s}, R\Gamma_c({\bf A}^1_{\overline{\bf F}_p}, {\cal
F}(Rf_{kn,!}s_{kn}^\ast {\cal L}_\psi)[kn-1]))^{-1}}{{\rm
det}(1-Fp^{-s}, ({\cal F}(Rf_{kn,!}s_{kn}^\ast {\cal
L}_\psi)[kn-1])|_{\bar 0})^{-1}}\\ &=&\frac{{\rm det}(1-F p^{-s},
R\Gamma_c(X_{k,\overline{\bf F}_p},s_{kn}^\ast {\cal
L}_\psi)(-1)[kn-2])^{-1}}{{\rm det}(1-Fp^{-s}, R\Gamma_c({\bf
G}^{kn}_{m,\overline{\bf F}_p}, s_{kn}^\ast {\cal
L}_\psi)[kn])^{-1}}.
\end{eqnarray*}
Let $Z_k$ be the subscheme of ${\bf P}^{kn-1}$ over ${\bf Z}$
defined by the conditions
$$x_{ij}\not=0,\; \sum_{j=1}^k \frac{1}{\prod_{i=1}^n x_{ij}}=0,$$
and let $Z_{k0}$ be the subscheme defined by
$$x_{ij}\not=0,\; \sum_{j=1}^k \frac{1}{\prod_{i=1}^n x_{ij}}=0,\;\sum_{i,j}
x_{ij}=0.$$ By Lemma 2.1, we have
\begin{eqnarray*}
{\rm det}(1-F p^{-s}, R\Gamma_c(X_{k,\overline{\bf F}_p},s_{kn}^\ast
{\cal L}_\psi)(-1)[kn-2])^{-1}&=&\frac{{\rm det}(1-F p^{-s},
R\Gamma_c(Z_{k0,\overline{\bf F}_p},\overline{\bf
Q}_l)(-2)[kn-4])^{-1}}{{\rm det}(1-F p^{-s},
R\Gamma_c(Z_{k,\overline{\bf
F}_p},\overline{\bf Q}_l)(-1)[kn-2])^{-1}}\\
&=&\frac{\zeta_{Z_{k0,{\bf F}_p}}(s-2)^{(-1)^{kn}}}{\zeta_{Z_{k,{\bf
F}_p}}(s-1)^{(-1)^{kn}}}.
\end{eqnarray*}
Let $Y_k$ be the subscheme of ${\bf P}^{kn-1}$ over ${\bf Z}$
defined by the condition
$$x_{ij}\not=0,$$ and let $Y_{k0}$ be the subscheme defined by
$$x_{ij}\not=0,\; \sum_{i,j}
x_{ij}=0.$$ By Lemma 2.1 again, we have
\begin{eqnarray*}
{\rm det}(1-Fp^{-s}, R\Gamma_c({\bf G}^{kn}_{m,\overline{\bf F}_p},
s_{kn}^\ast {\cal L}_\psi)[kn])^{-1} &=& \frac{{\rm det}(1-Fp^{-s},
R\Gamma_c(Y_{k0,\overline{\bf F}_p},\overline{\bf
Q}_l)(-1)[kn-2])^{-1}}{{\rm det}(1-Fp^{-s},
R\Gamma_c(Y_{k,\overline{\bf F}_p},\overline{\bf Q}_l) [kn])^{-1}}
\\
&=&\frac{\zeta_{Y_{k0,{\bf F}_p}}(s-1)^{(-1)^{kn}}}{\zeta_{Y_{k,{\bf
F}_p}}(s)^{(-1)^{kn}}}.
\end{eqnarray*}
So we finally get
\begin{eqnarray*}
L_p({\bf G}_{m,{\bf F}_p}, \otimes^k{\rm Kl}_{n+1}, s)&=& \frac{{\rm
det}(1-F p^{-s}, R\Gamma_c(X_{k,\overline{\bf F}_p},s_{kn}^\ast
{\cal L}_\psi)(-1)[kn-2])^{-1}}{{\rm det}(1-Fp^{-s}, R\Gamma_c({\bf
G}^{kn}_{m,\overline{\bf F}_p},
s_{kn}^\ast {\cal L}_\psi)[kn])^{-1}}\\
&=&\left(\frac{\zeta_{Z_{k0,{\bf F}_p}}(s-2)\zeta_{Y_{k,{\bf
F}_p}}(s)}{\zeta_{Z_{k,{\bf F}_p}}(s-1)\zeta_{Y_{k0,{\bf
F}_p}}(s-1)}\right)^{(-1)^{kn}}.
\end{eqnarray*}
Hence
$$\prod_p L_p({\bf G}_{m,{\bf F}_p}, \otimes^k{\rm Kl}_{n+1}, s)=
\left(\frac{\zeta_{Z_{k0}}(s-2)
\zeta_{Y_k}(s)}{\zeta_{Z_{k}}(s-1)\zeta_{Y_{k}}(s-1)}\right)^{(-1)^{kn}}.$$
This proves the assertions about the $L$-functions of $\otimes^k
{\rm Kl}_{n+1}$ in Theorem 0.1.

\bigskip
\bigskip
\centerline {\bf 3. The $L$-function of ${\rm Sym}^k{\rm Kl}_{n+1}$}
\bigskip
\bigskip
\noindent {\bf Lemma 3.1.} Let $V$ be a $\overline {\bf Q}_l$-vector
space, let $\pi:V\to V$ and $F:V\to V$ be two linear maps such that
$\pi^2=\pi$ and $F\pi=\pi F$. Then we have
$${\rm det}(1-Ft, {\rm im}(\pi))={\rm det}(1-F\pi t, V).$$

\bigskip
\noindent {\bf Proof.} Since $\pi^2=\pi$, we have
$$V={\rm ker}(\pi)\oplus {\rm im}(\pi),$$
and
$$\pi|_{{\rm ker}(\pi)}=0, \; \pi|_{{\rm im}(\pi)}={\rm id}.$$
Since $F\pi=\pi F$, the subspaces ${\rm ker}(\pi)$ and ${\rm
im}(\pi)$ are stable under $F$. It follows that
\begin{eqnarray*}
{\rm det}(1-F\pi t, V)
&=& {\rm det}(1-F\pi t,{\rm im}(\pi)){\rm det}(1-F\pi t, {\rm ker}(\pi))\\
&=&{\rm det}(1-Ft, {\rm im}(\pi)).
\end{eqnarray*}

\bigskip
Denote the coordinates of ${\bf G}_m^{kn}$ by $x_{ij}$ $(i=1,\ldots,
n, \; j=1,\ldots, k)$. Let
\begin{eqnarray*}
s_{kn}:{\bf G}_m^{kn}&\to&{\bf A}^1,\\
f_{kn}:{\bf G}_m^{kn}&\to&{\bf A}^1
\end{eqnarray*}
be the morphisms
\begin{eqnarray*}
s_{kn}((x_{ij}))=\sum_{i,j}x_{ij},\\
f_{kn}((x_{ij}))=\sum_{j=1}^k \frac{1}{\prod_{i=1}^n x_{ij}},
\end{eqnarray*}
respectively. Recall that in the previous section, we obtain the
isomorphisms (1) and (3):
$$
\otimes^k{\rm Kl}_{n+1}\cong \biggl({\cal F}(\ast^ k (j_! i^\ast{\rm
Kl}_n))[k-1]\biggr)|_{{\bf G}_m}\cong \biggl({\cal
F}(Rf_{kn,!}s_{kn}^\ast {\cal L}_\psi)[kn-1]\biggr)|_{{\bf G}_m}.$$
The group $\hbox{\frak S}_k$ acts on $\otimes^k {\rm Kl}_{n+1}$ and
on $\ast^ k (j_! i^\ast{\rm Kl}_n)$ by permuting the factors, and it
acts on $Rf_{kn,!}s_{kn}^\ast {\cal L}_\psi$ by permuting the
coordinates $x_{i1},\ldots, x_{ik}$ of ${\bf G}_m^{kn}$ for each
$i$. These actions are compatible with the above isomorphisms. By
Grothendieck's formula for $L$-functions, we have
$$L_p({\bf G}_{m,{\bf F}_p}, {\rm Sym}^k{\rm Kl}_{n+1}, s)
={\rm det}(1-Fp^{-s}, R\Gamma_c ({\bf G}_{m,\overline{\bf F}_p},
{\rm Sym}^k {\rm Kl}_{n+1}))^{-1}.$$  Let
$$\pi=\frac{1}{k_!}\sum_{\sigma\in\hbox{\frak S}_k} \sigma.$$
We have $\pi^2=\pi$, and $\pi$ induces the projection of $\otimes ^k
{\rm Kl}_{n+1}$ to its direct factor ${\rm Sym}^k{\rm Kl}_{n+1}$. It
follows that
\begin{eqnarray*}
H_c^m ({\bf G}_{m,\overline{\bf F}_p}, {\rm Sym}^k {\rm Kl}_{n+1})
&\cong &{\rm im}(H_c^m ({\bf G}_{m,\overline{\bf F}_p}, \otimes^k
{\rm Kl}_{n+1})\stackrel\pi\to H_c^m({\bf G}_{m,\overline{\bf F}_p},
\otimes^k {\rm Kl}_{n+1}))
\end{eqnarray*}
for all $m$. Applying Lemma 3.1 to $$\pi:H_c^m ({\bf
G}_{m,\overline{\bf F}_p}, \otimes^k{\rm Kl}_{n+1})\to H_c^m ({\bf
G}_{m,\overline{\bf F}_p}, \otimes^k{\rm Kl}_{n+1}),$$ we get
$${\rm det}(1-Fp^{-s}, H_c^m
({\bf G}_{m,\overline{\bf F}_p}, {\rm Sym}^k{\rm Kl}_{n+1}))={\rm
det}(1- F\pi p^{-s}, H_c^m ({\bf G}_{m,\overline{\bf F}_p},
\otimes^k {\rm Kl}_{n+1})).$$ It follows that
\begin{eqnarray*}
L_p({\bf G}_{m,{\bf F}_p}, {\rm Sym}^k{\rm Kl}_{n+1},s)&=& {\rm
det}(1- F\pi p^{-s}, R\Gamma_c ({\bf G}_{m,\overline{\bf
F}_p}, \otimes^k{\rm Kl}_{n+1}))^{-1}\\
&=&{\rm det}(1-F\pi p^{-s}, R\Gamma_c ({\bf G}_{m,\overline{\bf
F}_p}, ({\cal F}(Rf_{kn,!}s_{kn}^\ast {\cal L}_\psi)[kn-1])|_{{\bf
G}_{m,\overline{\bf
F}_p}}))^{-1}\\
&=& \frac{{\rm det}(1-F\pi p^{-s}, R\Gamma_c ({\bf
A}^1_{\overline{\bf F}_p}, {\cal F}(Rf_{kn,!}s_{kn}^\ast {\cal
L}_\psi)[kn-1]))^{-1}}{{\rm det}(1-F\pi p^{-s}, ({\cal
F}(Rf_{kn,!}s_{kn}^\ast {\cal L}_\psi)[kn-1])|_{\bar 0})^{-1}}.
\end{eqnarray*}
The same argument as in \S 2 shows that
\begin{eqnarray*}
&&{\rm det}(1-F\pi p^{-s}, R\Gamma_c({\bf A}^1_{\overline{\bf F}_p},
{\cal F}(Rf_{kn,!}s_{kn}^\ast {\cal L}_\psi)[kn-1]))^{-1}\\ &=&{\rm
det}(1-F\pi p^{-s}, R\Gamma_c(X_{k,\overline{\bf F}_p},s_{kn}^\ast
{\cal L}_\psi)(-1)[kn-2])^{-1}\\ &=&\frac{{\rm det}(1-F\pi p^{-s},
R\Gamma_c(Z_{k0,\overline{\bf F}_p},\overline{\bf
Q}_l)(-2)[kn-4])^{-1}}{{\rm det}(1-F\pi p^{-s},
R\Gamma_c(Z_{k,\overline{\bf F}_p},\overline{\bf
Q}_l)(-1)[kn-2])^{-1}},
\end{eqnarray*}
where $Z_k$ is the subscheme of ${\bf P}^{kn-1}$ over ${\bf Z}$
defined by the condition
$$x_{ij}\not=0,\; \sum_{j=1}^k \frac{1}{\prod_{i=1}^n x_{ij}}=0,$$
$Z_{k0}$ is the subscheme defined by
$$x_{ij}\not=0,\; \sum_{j=1}^k \frac{1}{\prod_{i=1}^n x_{ij}}=0,\;\sum_{i,j}
x_{ij}=0,$$ and the group $\hbox{\frak S}_k$ acts on $Z_k$ and on
$Z_{k0}$ by permuting the homogeneous coordinates $x_{i1},\ldots,
x_{ik}$ for each $i$. The same argument as in \S 2 also shows that
\begin{eqnarray*}
&&{\rm det}(1-F\pi p^{-s}, ({\cal F}(Rf_{kn,!}s_{kn}^\ast {\cal
L}_\psi)[kn-1])|_{\bar 0})^{-1}\\&=&{\rm det}(1-F\pi p^{-s},
R\Gamma_c({\bf G}^{kn}_{m,\overline{\bf F}_p}, s_{kn}^\ast {\cal
L}_\psi)[kn])^{-1}\\
&=& \frac{{\rm det}(1-F\pi p^{-s}, R\Gamma_c(Y_{k0,\overline{\bf
F}_p},\overline{\bf Q}_l)(-1)[kn-2])^{-1}}{{\rm det}(1-F\pi p^{-s},
R\Gamma_c(Y_{k,\overline{\bf F}_p},\overline{\bf Q}_l) [kn])^{-1}},
\end{eqnarray*}
where $Y_k$ is the subscheme of ${\bf P}^{kn-1}$ defined by
$$x_{ij}\not=0,$$  $Y_{k0}$ is the subscheme defined by
$$x_{ij}\not=0,\;\sum_{i,j}
x_{ij}=0,$$ and the group $\hbox{\frak S}_k$ acts on $Y_k$ and on
$Y_{k0}$ by permuting the homogeneous coordinates $x_{i1},\ldots,
x_{ik}$ for each $i$. So we have
\begin{eqnarray*}
&&L_p({\bf G}_{m,{\bf F}_p}, {\rm Sym}^k{\rm Kl}_{n+1},s)\\
&=& \frac{{\rm det}(1-F\pi p^{-s}, R\Gamma_c ({\bf
A}^1_{\overline{\bf F}_p}, {\cal F}(Rf_{kn,!}s_{kn}^\ast {\cal
L}_\psi)[kn-1]))^{-1}}{{\rm det}(1-F\pi p^{-s}, ({\cal
F}(Rf_{kn,!}s_{kn}^\ast {\cal L}_\psi)[kn-1])|_{\bar 0})^{-1}}\\
&=&\frac{{\rm det}(1-F\pi p^{-s}, R\Gamma_c(Z_{k0,\overline{\bf
F}_p},\overline{\bf Q}_l)(-2)[kn-4])^{-1}{\rm det}(1-F\pi p^{-s},
R\Gamma_c(Y_{k,\overline{\bf F}_p},\overline{\bf Q}_l)
[kn])^{-1}}{{\rm det}(1-F\pi p^{-s}, R\Gamma_c(Z_{k,\overline{\bf
F}_p},\overline{\bf Q}_l)(-1)[kn-2])^{-1}{\rm det}(1-F\pi p^{-s},
R\Gamma_c(Y_{k0,\overline{\bf F}_p},\overline{\bf
Q}_l)(-1)[kn-2])^{-1}}.
\end{eqnarray*}
Let $$a:Z_{k0}\to{\rm Spec}\,{\bf Z},\; b:Z_{k}\to{\rm Spec}\,{\bf
Z},\; c:Y_{k0}\to{\rm Spec}\,{\bf Z},\; d:Y_{k}\to{\rm Spec}\,{\bf
Z}$$  be the structure morphisms of $Z_{k0}$, $Z_k$, $Y_{k0}$ and
$Y_k$, respectively. By Lemma 3.1, we have
\begin{eqnarray*}
{\rm det}(1-F\pi p^{-s}, H^m_c(Z_{k0,\overline{\bf
F}_p},\overline{\bf Q}_l)(-2))&=&{\rm det}(1-F p^{-(s-2)}, {\rm
im}(H^m_c(Z_{k0,\overline{\bf F}_p},\overline{\bf
Q}_l)\stackrel{\pi}\to
H^m_c(Z_{k0,\overline{\bf F}_p},\overline{\bf Q}_l)))\\
&=& {\rm det}(1-F_p p^{-(s-2)}, {\rm im}(R^m a_! \overline{\bf
Q}_l\stackrel{\pi}\to R^ma_!\overline{\bf
Q}_l)) \\
{\rm det}(1-F\pi p^{-s}, H^m_c(Z_{k,\overline{\bf
F}_p},\overline{\bf Q}_l)(-1))&=&{\rm det}(1-F p^{-(s-1)}, {\rm im}(
H^m_c(Z_{k,\overline{\bf F}_p},\overline{\bf Q}_l)\stackrel{\pi}\to
H^m_c(Z_{k,\overline{\bf F}_p},\overline{\bf Q}_l)))\\
&=&{\rm det}(1-F_p p^{-(s-1)}, {\rm im}(R^mb_!\overline{\bf
Q}_l\stackrel{\pi}\to R^mb_!\overline{\bf Q}_l)),\\
{\rm det}(1-F\pi p^{-s}, H^m_c(Y_{k0,\overline{\bf
F}_p},\overline{\bf Q}_l)(-1)) &=&{\rm det}(1-F p^{-(s-1)}, {\rm
im}(H^m_c(Y_{k0,\overline{\bf F}_p},\overline{\bf
Q}_l)\stackrel{\pi}\to
H^m_c(Y_{k0,\overline{\bf F}_p},\overline{\bf Q}_l)))\\
&=& {\rm det}(1-F_p p^{-(s-1)}, {\rm im}(R^mc_!\overline{\bf
Q}_l\stackrel{\pi}\to
R^mc_!\overline{\bf Q}_l)),\\
{\rm det}(1-F\pi p^{-s}, H^m_c(Y_{k,\overline{\bf
F}_p},\overline{\bf Q}_l)) &=&{\rm det}(1-F p^{-s}, {\rm
im}(H^m_c(Y_{k,\overline{\bf F}_p},\overline{\bf
Q}_l)\stackrel{\pi}\to
H^m_c(Y_{k,\overline{\bf F}_p},\overline{\bf Q}_l)))\\
&=&{\rm det}(1-F_p p^{-s},{\rm im}(R^md_!\overline{\bf
Q}_l\stackrel{\pi}\to R^md_!\overline{\bf Q}_l)).
\end{eqnarray*}
So we have
\begin{eqnarray*}
&& L_p({\bf G}_{m,{\bf F}_p}, {\rm Sym}^k{\rm Kl}_{n+1}, s)\\
&=&\frac{{\rm det}(1-F\pi p^{-s}, R\Gamma_c(Z_{k0,\overline{\bf
F}_p},\overline{\bf Q}_l)(-2)[kn-4])^{-1}{\rm det}(1-F\pi p^{-s},
R\Gamma_c(Y_{k,\overline{\bf F}_p},\overline{\bf Q}_l)
[kn])^{-1}}{{\rm det}(1-F\pi p^{-s}, R\Gamma_c(Z_{k,\overline{\bf
F}_p},\overline{\bf Q}_l)(-1)[kn-2])^{-1}{\rm det}(1-F\pi p^{-s},
R\Gamma_c(Y_{k0,\overline{\bf F}_p},\overline{\bf
Q}_l)(-1)[kn-2])^{-1}}\\
&=&\prod_m\left(\frac{{\rm det}(1-F_p p^{-(s-2)}, {\rm im}(R^m a_!
\overline{\bf Q}_l\stackrel{\pi}\to R^ma_!\overline{\bf Q}_l)){\rm
det}(1-F_p p^{-s},{\rm im}(R^md_!\overline{\bf Q}_l\stackrel{\pi}\to
R^md_!\overline{\bf Q}_l))}{{\rm det}(1-F_p p^{-(s-1)}, {\rm
im}(R^mb_!\overline{\bf Q}_l\stackrel{\pi}\to R^mb_!\overline{\bf
Q}_l)){\rm det}(1-F_p p^{-(s-1)}, {\rm im}(R^mc_!\overline{\bf
Q}_l\stackrel{\pi}\to R^mc_!\overline{\bf
Q}_l))}\right)^{(-1)^{kn+m+1}},
\end{eqnarray*}
and
\begin{eqnarray*}
&& \prod_p L_p({\bf G}_{m,{\bf F}_p}, {\rm Sym}^k{\rm Kl}_{n+1},
s)\\
&=& \prod_m\left(\frac{L({\rm Spec}\, {\bf Z}, {\rm im}(R^ma_!
\overline {\bf Q}_l\stackrel{\pi}\to R^m a_!\overline {\bf Q}_l),
s-2)L({\rm Spec}\, {\bf Z}, {\rm im}(R^md_! \overline {\bf
Q}_l\stackrel{\pi}\to R^m d_!\overline {\bf Q}_l), s)}{L({\rm
Spec}\, {\bf Z}, {\rm im}(R^mb_! \overline {\bf
Q}_l\stackrel{\pi}\to R^mb_! \overline {\bf Q}_l ), s-1)L({\rm
Spec}\, {\bf Z}, {\rm im}(R^mc_! \overline {\bf
Q}_l\stackrel{\pi}\to R^m c_!\overline {\bf Q}_l), s-1)}
\right)^{(-1)^{kn+m}}.
\end{eqnarray*}
The above sheaf ${\rm im}(R^ma_! \overline {\bf Q}_l\stackrel\pi\to
R^m a_!\overline {\bf Q}_l)$ and the similar sheaves for the
morphisms $b$, $c$ and $d$ can be made more explicit. The group
$\hbox{\frak S}_k$ acts on $R^ma_!\overline{\bf Q}_l$. We have
$$(R^ma_!\overline {\bf Q}_l)^{\hbox {\frak S}_k}\cong
{\rm im}(R^ma_! \overline {\bf Q}_l\stackrel\pi\to R^m a_!\overline
{\bf Q}_l).$$ Let $a':Z_{k0}/\hbox{\frak S}_k\to {\rm Spec}\, {\bf
Z}$ be the structure morphism of the quotient of $Z_{k0}$ by $\hbox
{\frak S}_k$. Then we have
$$(R^ma_!\overline {\bf Q}_l)^{\hbox {\frak S}_k}\cong R^m
a'_!\overline {\bf Q}_l.$$ To prove this, we use the
Hochschild-Serre type spectral sequences in [G] 5.2.1. These
spectral sequences are constructed by Grothendieck for the
cohomology of sheaves of abelian groups on topological spaces. We
can construct similar spectral sequences for the cohomology of
\'etale sheaves of torsion abelian groups on schemes. We then use
the fact that $H^i(\hbox{\frak S}_k, -)$ are annihilated by $k!$ for
all $i>0$ to conclude that similar spectral sequences degenerate for
cohomology of $\overline {\bf Q}_l$-sheaves. So we have
$$R^m
a'_!\overline {\bf Q}_l\cong {\rm im}(R^ma_! \overline {\bf
Q}_l\stackrel\pi\to R^m a_!\overline {\bf Q}_l).$$ Therefore we have
\begin{eqnarray*}
&& L_p({\bf G}_{m,{\bf F}_p}, {\rm Sym}^k{\rm Kl}_{n+1}, s)\\
&=&\prod_m\left(\frac{{\rm det}(1-F_p p^{-(s-2)}, {\rm im}(R^m a_!
\overline{\bf Q}_l\stackrel{\pi}\to R^ma_!\overline{\bf Q}_l)){\rm
det}(1-F_p p^{-s},{\rm im}(R^md_!\overline{\bf Q}_l\stackrel{\pi}\to
R^md_!\overline{\bf Q}_l))}{{\rm det}(1-F_p p^{-(s-1)}, {\rm
im}(R^mb_!\overline{\bf Q}_l\stackrel{\pi}\to R^mb_!\overline{\bf
Q}_l)){\rm det}(1-F_p p^{-(s-1)}, {\rm im}(R^mc_!\overline{\bf
Q}_l\stackrel{\pi}\to R^mc_!\overline{\bf
Q}_l))}\right)^{(-1)^{kn+m+1}}\\
&=&\prod_m\left(\frac{{\rm det}(1-F_p p^{-(s-2)}, R^m
a'_!\overline{\bf Q}_l){\rm det}(1-F_p p^{-s},R^md'_!\overline{\bf
Q}_l)}{{\rm det}(1-F_p p^{-(s-1)}, R^mb'_!\overline{\bf Q}_l){\rm
det}(1-F_p p^{-(s-1)},
R^mc'_!\overline{\bf Q}_l)}\right)^{(-1)^{kn+m+1}}\\
&=&\left(\frac{\zeta_{Z_{k0,{\bf F}_p}/\hbox{\frak
S}_k}(s-2)\zeta_{Y_{k,{\bf F}_p}/\hbox{\frak
S}_k}(s)}{\zeta_{Z_{k,{\bf F}_p}/\hbox{\frak
S}_k}(s-1)\zeta_{Y_{k0,{\bf F}_p}/\hbox{\frak
S}_k}(s-1)}\right)^{(-1)^{kn}}.
\end{eqnarray*}
This proves the assertions about the $L$-functions of ${\rm
Sym}^k{\rm Kl}_{n+1}$ in Theorem 0.2.

\bigskip
Similarly, by working with
$$\pi'=\frac{1}{k_!}\sum_{\sigma\in\hbox{\frak S}_k}{\rm sgn}(\sigma)\sigma$$
instead of $\pi$, we can prove the following result for the $k$-th
exterior power.

\bigskip
\noindent {\bf Theorem 3.2.} Notation as above.  We have
\begin{eqnarray*}
&& L_p({\bf G}_{m,{\bf F}_p}, \bigwedge^k{\rm Kl}_{n+1}, s)\\
&=&\prod_m\left(\frac{{\rm det}(1-F_p p^{-(s-2)}, {\rm im}(R^m a_!
\overline{\bf Q}_l\stackrel{\pi'}\to R^ma_!\overline{\bf Q}_l)){\rm
det}(1-F_p p^{-s},{\rm im}(R^md_!\overline{\bf
Q}_l\stackrel{\pi'}\to R^md_!\overline{\bf Q}_l))}{{\rm det}(1-F_p
p^{-(s-1)}, {\rm im}(R^mb_!\overline{\bf Q}_l\stackrel{\pi'}\to
R^mb_!\overline{\bf Q}_l)){\rm det}(1-F_p p^{-(s-1)}, {\rm
im}(R^mc_!\overline{\bf Q}_l\stackrel{\pi'}\to R^mc_!\overline{\bf
Q}_l))}\right)^{(-1)^{kn+m+1}},
\end{eqnarray*}
and
\begin{eqnarray*}
&& \prod_p L_p({\bf G}_{m,{\bf F}_p}, \bigwedge^k{\rm Kl}_{n+1},
s)\\
&=& \prod_m\left(\frac{ L({\rm Spec}\, {\bf Z}, {\rm im}(R^ma_!
\overline {\bf Q}_l\stackrel{\pi'}\to R^m a_!\overline {\bf Q}_l),
s-2)L({\rm Spec}\, {\bf Z}, {\rm im}(R^md_! \overline {\bf
Q}_l\stackrel{\pi'}\to R^m d_!\overline {\bf Q}_l), s)}{L({\rm
Spec}\, {\bf Z}, {\rm im}(R^mb_! \overline {\bf
Q}_l\stackrel{\pi'}\to R^mb_! \overline {\bf Q}_l ), s-1)L({\rm
Spec}\, {\bf Z}, {\rm im}(R^mc_! \overline {\bf
Q}_l\stackrel{\pi'}\to R^m c_!\overline {\bf Q}_l), s-1)}
\right)^{(-1)^{kn+m}}.
\end{eqnarray*}

\bigskip
The sheaf ${\rm im}(R^ma_! \overline {\bf Q}_l\stackrel{\pi'}\to R^m
a_!\overline {\bf Q}_l)$ and the similar sheaves for the morphisms
$b$, $c$ and $d$ can again be made more explicit. Let ${\cal S}$ be
the constant sheaf $\overline {\bf Q}_l$ on $Z_{k0}/\hbox{\frak
S}_k$ provided with an action of $\hbox {\frak S}_k$ so that
$\sigma\in \hbox{\frak S}_k$ acts as multiplication by ${\rm
Sgn}(\sigma)$. Let $p_{Z_{k0}}:Z_{k0}\to Z_{k0}/\hbox{\frak S}_k$ be
the projection. Using Hochschild-Serre type spectral sequences, one
can show that
$$R^m
a'_!\left((p_{Z_{k0},\ast}\overline {\bf Q}_l\otimes {\cal
S})^{\hbox{\frak S}_k}\right)\cong {\rm im}(R^ma_! \overline {\bf
Q}_l\stackrel{\pi'}\to R^m a_!\overline {\bf Q}_l).$$

\bigskip
\centerline{\bf 4. Simplified formulas}
\bigskip

The formula for $\prod_p L_p({\bf G}_{m,{\bf F}_p}, \otimes^k{\rm
Kl}_{n+1}, s)$ in Theorem 0.2 can be significantly simplified. Since
$Y_k$ is isomorphic to ${\bf G}_m^{kn-1}$, we have $$\# Y_k({\bf
F}_q) =(q-1)^{kn-1}.$$ A simple inclusion-exclusion argument shows
that
$$\# Y_{k0}({\bf F}_q) =\frac{1}{q}\left( (q-1)^{kn-1} +(-1)^{kn}\right).$$
This gives the relation
$$\frac{\zeta_{Y_k}(s)}{\zeta_{Y_{k0}}(s-1)} =\zeta(s)^{(-1)^{kn}},$$
where $\zeta(s)$ is the Riemann zeta function. Similarly, one checks
that
$$\# Z_{k}({\bf F}_q) =(q-1)^{k(n-1)}\frac{1}{q}\left( (q-1)^{k-1} +(-1)^{k}\right)
=\frac{1}{q}\left( (q-1)^{kn-1} +(-1)^{k}(q-1)^{k(n-1)}\right).$$
Thus $\zeta_{Z_k}(s)$ is also determined explicitly by the Riemann
zeta function. The only non-trivial factor in the formula for
$\prod_p L_p({\bf G}_{m,{\bf F}_p}, \otimes^k{\rm Kl}_{n+1}, s)$ is
the zeta function $\zeta_{Z_{k0}}(s)$. From the last equation
defining $Z_{k0}$, we get
$$x_{nk}= -\left(\sum_{i=1}^{n-1}x_{ik}+\sum_{i=1}^n\sum_{j=1}^{k-1} x_{ij}\right).$$
Substituting this into the second equation defining $Z_{k0}$, we see
that $Z_{k0}$ is isomorphic to the toric hypersurface $W_k$ in
$$\{[x_{ij}]\in{\bf P}^{kn-1}|x_{11}=1,x_{ij}\not =0\}\cong {\bf G}_m^{kn-1}$$ defined by
$$x_{11}=1,\;\sum_{j=1}^{k-1}\frac{1}{\prod_{i=1}^n x_{ij}}
\left(\sum_{i=1}^{n-1}x_{ik}+\sum_{i=1}^n\sum_{j=1}^{k-1}
x_{ij}\right)-\frac{1}{\prod_{i=1}^{n-1}x_{ik}}=0.$$ Thus, we obtain
the simplified formula
$$\prod_p L_p({\bf G}_{m,{\bf F}_p},
\otimes^k{\rm Kl}_{n+1}, s) = \zeta(s) \left(
\frac{\zeta_{W_k}(s-2)}{\zeta_{Z_{k}}(s-1)}\right)^{(-1)^{kn}}.$$

The formula for the $L$-function of ${\rm Sym}^k {\rm Kl}_{n+1}$ is
more complicated. The scheme $Y_k/\hbox{\frak S}_k$ can be
explicitly described as follows. Let $S=k[x_{ij}]$ be the polynomial
ring with the canonical grading by the degrees of polynomials. The
group $\hbox{\frak S}_k$ acts on $S$ by permuting the indeterminates
$x_{i1},\ldots, x_{ik}$. Let $f=\prod_{i,j} x_{ij}$. Then $Y_k={\rm
Spec}\, S_{(f)}$. Let $s_{ij}$ be the $j$-th elementary symmetric
polynomial of $x_{i1},\ldots, x_{ik}$. Then the subring of $S$ fixed
by $\hbox {\frak S}_k$ is
$$S^{\hbox{\frak S}_k}=k[s_{ij}].$$
Let $S'=k[s_{ij}]$. It is isomorphic to a polynomial ring. Introduce
a grading on $S'$ by setting ${\rm deg}(s_{ij})=j$. Then we have
$$(S_{(f)})^{\hbox{\frak S}_k}=S'_{(f)}$$ and hence
$$Y_k/\hbox{\frak S}_k={\rm Spec}\, S'_{(f)}.$$
Let ${\bf Q}^{kn-1}={\rm Proj}\, S'$ which is a weighted projective
space. Then $Y_k/\hbox{\frak S}_k$ is the complement of the
hypersurface $f=0$ in ${\bf Q}^{kn-1}$.

\bigskip
\noindent {\bf References}

\bigskip

\bigskip
\noindent [CE] T. Choi and R. Evans, {\it Congruences for powers of
Kloosterman sums}, Intern. J. Number Theory, 3 (2007), 105-117.

\bigskip
\noindent [Ev] R. Evans, {\it Seventh power moments of Kloosterman
sums}, to appear.

\bigskip
\noindent [FW1] L. Fu and D. Wan, {\it L-functions for symmetric
products of Kloosterman sums}, J. Reine Angew. Math., 589 (2005),
79-103.

\bigskip
\noindent [FW2] L. Fu and D. Wan, {\it Trivial factors for
L-functions of symmetric products of Kloosterman sheaves}, Finite
Fields and Their Applications, to appear.

\bigskip
\noindent [FH] W. Fulton and J. Harris, {\it Representation theory,
a first course}, Springer-Verlag (1991).

\bigskip
\noindent [G] A. Grothendieck, {\it Sur quelques points d'alg\`ebre
homologique}, T\^ohoku Math. J. 9 (1957), 119-221.

\bigskip
\noindent [GK] E. Grosse-Kl\"onne, {\it On families of pure slope
$L$-functions}, Doc. Math.,  8 (2003), 1-42.

\bigskip
\noindent [HS] K. Hulek,  J. Spandaw, B. van Geemen,  D. van
Straten, {\it The modularity of the Barth-Nieto quintic and its
relatives}, Adv. Geom. 1 (2001), no. 3, 263-289.

\bigskip
\noindent [L] G. Laumon, {\it Transformation de Fourier, constantes
d'\'equations fontionnelles, et conjecture de Weil}, Publ. Math.
IHES 65 (1987), 131-210.

\bigskip
\noindent [PTV] C. Peters, J. Top and M. van der Vlugt, {\it The
Hasse zeta function of a K3 surface related to the number of words
of weight $5$ in the Melas codes}, J. Reine Angew. Math., 432
(1992), 151-176.

\bigskip
\noindent [Ro] P. Robba, {\it Symmetric powers of $p$-adic Bessel
equation}, J. Reine Angew. Math., 366 (1986), 194-220.

\bigskip
\noindent [SGA $4\frac{1}{2}$] P. Deligne et al, {\it Cohomologie
\'Etale}, Lecture Notes in Math., 569, Springer-Verlag (1977).

\bigskip
\noindent [W1] D. Wan, {\it Dwork's conjecture on unit root zeta
functions}, Ann. Math., 150 (1999), 867-927.

\end{document}